\documentclass[a4paper, 10pt]{article}

\usepackage{amsmath,amsfonts}

\usepackage[T1]{fontenc}

\newtheorem{Theorem}{\quad Theorem}[section]

\title {A uniform result in dimension 2.}

\author {Samy Skander Bahoura \footnote {E-mails: samybahoura@yahoo.fr, samybahoura@gmail.com}}

\date{ \small {Universit\'e Pierre et Marie Curie, 75005, Paris, France.}}

\begin{document}

\maketitle

\begin{abstract} 

We give a uniform result in dimension 2 for the solutions to an equation on compact Riemannian surface without boundary.

\end{abstract}

\smallskip


\section{Introduction  and Main result}

\bigskip

We set $ \Delta=-\nabla_i(\nabla^i) $ the Laplace-Beltrami operator. We are on compact Riemannian surface $ (M,g) $ without boundary.

\bigskip

We start with the following example: for all $ \epsilon >0 $  the constant functions $ z_{\epsilon}=\log\dfrac{\epsilon}{a} $ with $ a>0 $, are solutions to $ \Delta z_{\epsilon}+\epsilon=ae^{z_{\epsilon}} $  and tend to $ -\infty $  uniformly on $ M $.

\bigskip

Question: What's about the solutions $ u_{\epsilon} $ to the following equation

$$ \Delta u_{\epsilon}+\epsilon=V_{\epsilon}e^{u_{\epsilon}}, \qquad (E_{\epsilon})$$

with $ 0 < a\leq V_{\epsilon}(x)\leq b<+\infty $ on $ M $ ?

\bigskip

Next, we assume $ V_{\epsilon} $ H\"olderian and $ V_{\epsilon} \to V $ in $ L^{\infty} $ 

\bigskip

The equation $ (E_{\epsilon}) $  is of prescribed scalar curvature type equation. The term $ \epsilon $ replace the scalar curvature.

\bigskip

\begin{Theorem}. If $ \epsilon \to 0 $, the solutions $ u_{\epsilon} $ to $ (E_{\epsilon}) $ satisfy:

$$ \sup_M u_{\epsilon} \to -\infty . $$

\end{Theorem}

By using the same arguments of the next theorem, we have:

\begin{Theorem}. If $ \epsilon \to 0 $, the solutions $ u_{\epsilon} $ to $ (E_{\epsilon}) $ satisfy:

$$ u_{\epsilon}-\log \epsilon \to k \in {\mathbb R}. $$

uniformly on $ M $.

\end{Theorem}

Thus, we have a unifrom bound for the solutions:

$$ k_1+\log \epsilon \leq u_{\epsilon} \leq \log \epsilon +k_2. $$

We also have another proof of the uniqueness result which appear in [2]. This proof uses Brezis Merle arguments.

\bigskip

\begin{Theorem}.  If $ \epsilon \to 0 $, the solutions $ u_{\epsilon} $ to $ (E_{\epsilon}) $ with $ V_{\epsilon}\equiv 1 $, are such:

$$ u_{\epsilon}\equiv \log \epsilon. $$

\end{Theorem}

\bigskip

\section{Proof of the theorems 1,2,3.}

\bigskip

\underbar {Proof of theorem 1:}

\bigskip

We have:

$$ \int_M V_{\epsilon} e^{u_{\epsilon}} \to 0. \qquad (*) $$

\bigskip

Let's consider $ x_{\epsilon} $ a point such that $ \max_M u_{\epsilon}=u_{\epsilon}(x_{\epsilon}) $, then $ x_{\epsilon} \to x_0 $.

\bigskip

We consider a neighborhood of $ x_0 $  and we use isothermal coordinates around $ x_0 $ (see [6]), there exists $ \alpha >0 $  and a regular function $ \phi $ such that:

$$ \Delta_{\cal E} u_{\epsilon}+\epsilon e^{\phi}=V_{\epsilon} e^{\phi}e^{u_{\epsilon}} \,\,\, {\rm in} \,\, B(0,\alpha). $$

The metric $ g $ of $ M $ satisfies $ g=e^{\phi}(dx^2+dy^2) $.

\smallskip

Let's consider $ u_0 $ such that:

$$ \Delta_{\cal E} u_0=e^{\phi} \,\,\, {\rm in } \,\,\, B(0,\alpha) .$$

(with Dirichlet condition for example).

The function $ v_{\epsilon}=u_{\epsilon}+\epsilon u_0 $ satisfies:

$$ \Delta_{\cal E} v_{\epsilon}=\tilde V_{\epsilon} e^{v_{\epsilon}}, $$

with $ \tilde V_{\epsilon}=V_{\epsilon} e^{\phi-\epsilon u_0} $. We use $ (*) $ to have:

$$ \int_{B(0,\alpha)} e^{v_{\epsilon}} \to 0 \,\,\,{\rm and} \,\,\, 0 < \tilde a \leq \tilde V_{\epsilon} \leq \tilde b. \qquad (**) $$

with $ \alpha >0 $.

\bigskip

The sequence $ v_{\epsilon} $ satisfies all the conditions of the theorem of Brezis and Merle, see [4].

\bigskip

As $ v_{\epsilon} $ satisfy $ (**) $, the last condtion of the theorem of [4] is not possible.

\bigskip

Now, suppose that the first assertion of the theorem of Brezis and Merle is true. We have the local boundedness result. We can say that $ u_{\epsilon} $ converge uniformly on $ M $ to  a fonction $ u $ and in $ C^2 $ topology by the elliptic estimates.

\bigskip

If we tend $ \epsilon  $ to $ 0 $ we get that $ u $ satisfies in the sense of distributions:

$$ \Delta u=Ve^u. $$

\smallskip

If we integrate the equation, we have a contradiction (since $ 0< a \leq V \leq b <+\infty $.

\bigskip

Thus, $ u_{\epsilon} $ satisfies the second assertion of the theorem of [4] and thus $ u_{\epsilon} $ diverge uniformly to  $ -\infty $ on $ M $.

\bigskip

\underbar {Proof of Theoreme 2,3:}

We set,

$$ w_{\epsilon} = u_{\epsilon} -\log \epsilon. $$

Then, $ w_{\epsilon} $ is solution to:

$$ \Delta w_{\epsilon} +\epsilon = \epsilon V_{\epsilon} e^{w_{\epsilon}}. $$

We use Brezis and Merle's theorem and the previous arguments of theoerm 1, to have a convergence to a constant:

$$ w_{\epsilon} \to w_{\infty}=ct, $$

uniformly on $ M $.

In isothermal coordinates around  $ x_0=\lim x_{\epsilon} $ with $ x_{\epsilon} $ such that, $ w_{\epsilon}(x_{\epsilon})=\max_M w_{\epsilon} $, as in the previous case

$$ \Delta_{\cal E} w_{\epsilon}+\epsilon e^{\phi}=\epsilon e^{\phi} V_{\epsilon}e^{w_{\epsilon}} \,\,\, {\rm in} \,\, B(0,\alpha). $$

The metric $ g $ of $ M $ satisfies $ g=e^{\phi}(dx^2+dy^2) $.
Let's consider $ u_0 $ such that:

$$ \Delta_{\cal E} u_0=e^{\phi} \,\,\, {\rm in } \,\,\, B(0,\alpha) .$$

(with Dirichlet condition for example).

The function $ v_{\epsilon}=w_{\epsilon}+\epsilon u_0 $ satisfies:

$$ \Delta_{\cal E} v_{\epsilon}=\tilde V_{\epsilon} e^{v_{\epsilon}}, $$

with $ \tilde V_{\epsilon}=\epsilon V_{\epsilon} e^{\phi-\epsilon u_0} $. 

One can apply the theorem of Brezis and Merle, see [4].

First we have,

$$ \int_M V_{\epsilon} e^{w_{\epsilon}} =|M|, $$

which imply,

$$ w_{\epsilon} \not \to -\infty, $$

and,

$$ \int_M \epsilon V_{\epsilon} e^{w_{\epsilon}} =\epsilon |M| \to 0, $$

which imply the non-concentration.

And,

$$ w_{\epsilon} \to w_{\infty}, $$

in the $ C^2 $ topology with,

$$ \Delta w_{\infty} =0 \Rightarrow w_{\infty} \equiv k\in{\mathbb R}. $$

For the third theorem we have:

$$ \int_M e^{w_{\epsilon}}=|M| \Rightarrow k=0. $$

We write:

$$ w_{\epsilon}=\bar w_{\epsilon} + f_i, $$

with the fact that,

$$ \int_M f_i=0, $$

and, $ f_i $ is solution to:

$$ \Delta f_i =\epsilon_i (e^{\bar w_i +f_i} -1)=\epsilon_i (e^{\bar w_i}(1+f_i+O(f_i^2))-1), $$

We multiply the equation by $ f_i $ and we integrate, we obtain:

$$ ||\nabla f_i||_{L^2}^2=o(||f_i||_{L^2}^2). $$

This is in contradiction with the Poincar\'e inequality if $ f_i\not \equiv 0 $.

Thus,

$$ f_i\equiv 0,\,\, w_i\equiv \bar w_i =0. $$

\end{document}